# Non-stationary extremal eigenvalue approximations in iterative solutions of linear systems and estimators for relative error

Divya Anand Subba and Murugesan Venkatapathi*

*Supercomputer Education and Research Center, Indian Institute of Science, Bangalore – 560012, India*
* Corresponding author: murugesh@serc.iisc.ernet.in

**Abstract**

Non-stationary approximations of the final value of a converging sequence are discussed, and we show that extremal eigenvalues can be reasonably estimated from the CG iterates without much computation at all. We introduce estimators of relative error for conjugate gradient (CG) methods that adopt past work on computationally efficient bounds of the absolute errors using quadrature formulas. The evaluation of the Gauss quadrature based estimates though, depends on a priori knowledge of extremal eigenvalues; and the upper bounds in particular that are useful as a stopping criterion fail in the absence of a reasonable *underestimate* of smallest eigenvalue. Estimators for relative errors in $A$-norm and their extension to errors in $l_2$ norm are presented with numerical results. Estimating the relative error from the residue in an iterative solution is required for efficient solution of a large problem with even a moderately high condition. Specifically, in a problem of solving for vector $x$ in $Ax=b$, the uncertainty between the strict upper bound in relative error [$\kappa\times||r^i||/||b||$] and its strict lower bound [$||r^i||/(\kappa\times||b||)$] is a factor of $\kappa^2$ (given residue $r^i= b-Ax^i$ is the residual vector at $i^{th}$ iteration and $\kappa$ the condition number of the square matrix $A$).

## 1. Introduction

In numerical models of linear systems, the stopping criterion of an iterative solution is typically based on the relative residue at an iteration; for the problem of solving for vector $x$ in $Ax=b$, this is given by $||r^i||/||b||$ where $r^i= b-Ax^i$ is the residual vector at $i^{th}$ iteration. But the error in $x$ is the indicator of the convergence behavior and the accuracy of the solution obtained; the relative error in $x$ preferred typically when the unknown $||x||$ varies significantly by the problem. Specifically for a square matrix $A$, it can be easily shown that the relative error is bounded such that $\kappa_A \times ||r^i||/||b|| \geq ||x-x^i||/||x|| \geq ||r^i||/(\kappa_A \times ||b||)$ where $\kappa_A$ is the condition number of $A$. If one defines a condition of the specific problem $\kappa_{(A,b)}$ (i.e.) given a matrix $A$ and the vector $b$, then $||x-x^i||/||x|| \rightarrow \kappa_{(A,b)} \times ||r^i||/||b||$ as $x^i \rightarrow x$. Computation of the condition numbers $\kappa_A$ or $\kappa_{(A,b)}$ are arithmetic operations of $O(N^3)$ for a general matrix and thus prohibitive when the dimension of matrix $N$ is large [1-3]; also the above loose bounds of the relative error are not useful as an efficient stopping criterion even if the condition of the matrix $\kappa_A$ is known, but $\kappa_A >> 1$. As an example, figure (1) shows the relative residue and the actual relative errors for two different vectors $b$ and a symmetric positive definite (SPD) matrix $A$ using the conjugate gradient algorithm. In one case, the relative residue is a few orders of magnitude smaller than the relative error, and alternately larger in the other case. A stopping criteria based on the residue alone can thus be optimistic and erroneous, or conservative and inefficient, depending on the condition of the specific problem $\kappa_{(A,b)}$, whether pre-conditioned or otherwise. The larger the problem and its condition number, more useful is the estimate of the relative error in enforcing a stopping criterion that is both cost effective but satisfying the required accuracy of solution. Here we extend the work on estimates of absolute error in CG iterates using quadrature rules [4, 5] that demand negligible computation, to include relative errors, specifically to the realistic scenario where the extremal eigenvalues are unknown. An efficient non-stationary estimation of the extremal eigenvalues in $O(N)$ operations at every iteration of the CG algorithm is suggested; here we use the iterative algorithm solving the problem to efficiently evaluate extremal eigenvalues as well. Other methods of estimating extremal eigenvalues like probabilistic bounds using a Krylov subspace of vector $b$ [6] are either dependent on vector $b$ or require arithmetic operations $>O(N^2)$ on their own. Also there exist anti-Gauss quadrature based methods for estimation of errors that do not need the smallest eigenvalue, but negligible computation and upper bounds are not guaranteed [7, 8]. The property of $A$-orthogonality of the CG increments is used to extend absolute error estimates to relative errors and numerical estimates of the relative error in $A$-norm are shown. The above method is extended to estimating the relative errors in $l_2$ norm; this is significant for applications where the $A$-norm/energy-norm is not an appropriate criterion for stopping.

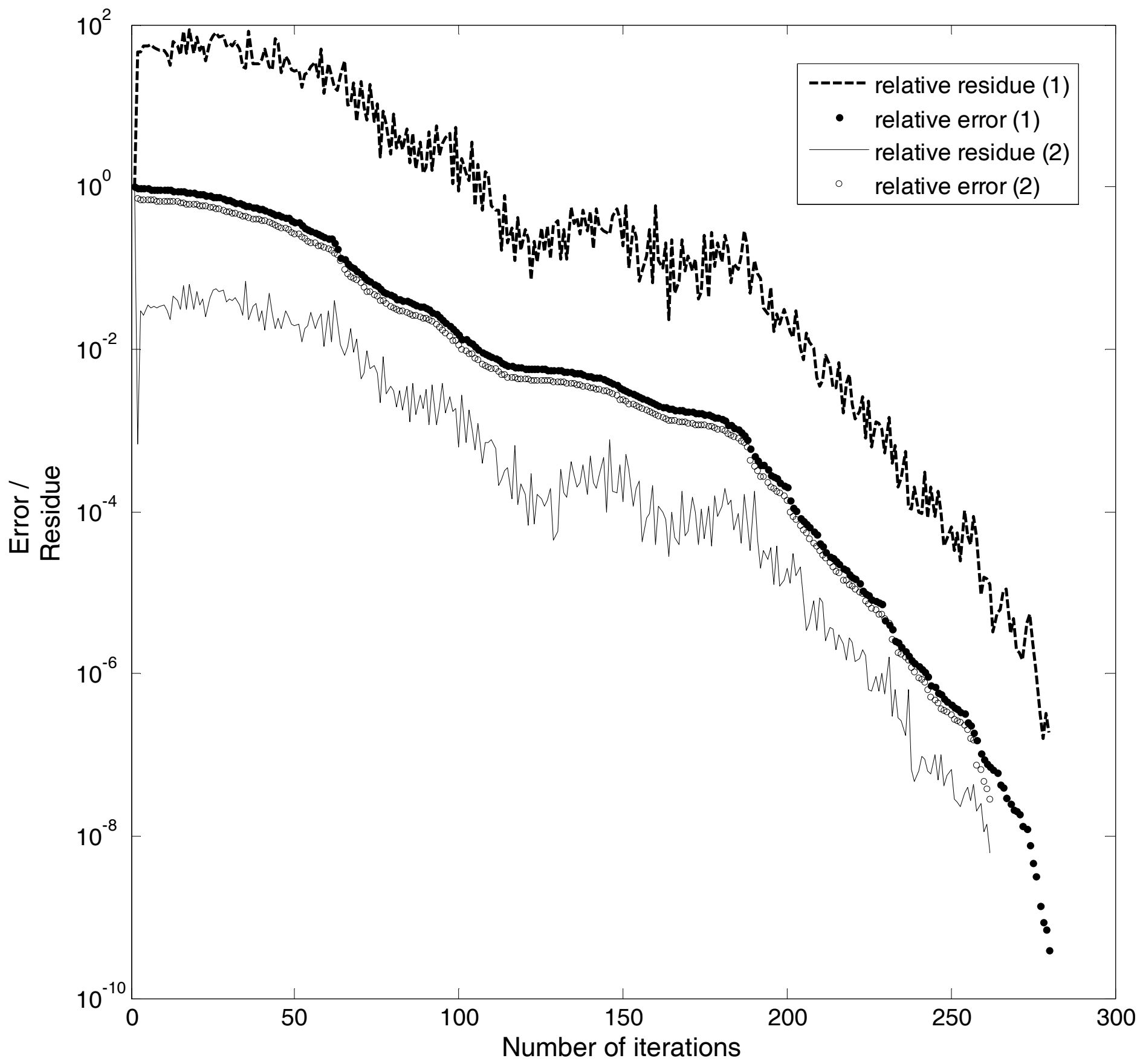

Figure 1: Relative errors ($||x\text{-}x^i||_A/||x||_A$) and residues ($||r^i||_A/||b||_A$) for two different vectors $b$ and a SPD matrix $A$ (bcsstk05; math.nist.gov/MatrixMarket). Case 1: vector $b$ is a sum of canonical basis (all 1s) for generality and Case 2: vector $b$ is eigenvector of smallest eigenvalue with a small (< 0.01) random perturbation of its values; $\kappa_A = 14281$, $N=153$.

## 2. Mathematical preliminaries

To use quadrature formulas in estimation of the absolute error, the quadratic form of the square of error in $A$-norm, $||x\text{-}x_k||_A^2 = r^TA^{-1}r$, can be represented by a Riemann–Stieltjes integral as below.

when $AA^T = A^T A$

$$\|\varepsilon\|_A^2 = I[A,\, r] \;=\; r^T A^{-1} r \;=\; r^T Q \Lambda^{-1} Q^T r$$

$$= q^T \Lambda^{-1} q = \sum_{i=1}^{n} \lambda_i^{-1} q_i^{\,2} = \int_a^b \lambda^{-1} d\alpha(\lambda)$$ where α is piecewise costant and defined

$$\alpha(\lambda) = 0 \;\; if \;\; \lambda \le a = \lambda \qquad \text{Eq. (1)}$$

$$= \sum_{j=1}^{i} q_j^{\,2} \;\; if \;\; \lambda_i \le \lambda < \lambda_{i+1}$$

$$= \sum_{j=1}^{n} q_j^{\,2} \;\; if \;\; \lambda \ge b = \lambda_n$$

Here $a$, $b$ (unlike vector $b$) are the smallest and largest eigenvalues of $A$ respectively. This allows us to use the Gauss, Gauss-Radau and Gauss-Lobatto formulas for a function $f$ given by

$$\int_a^b f(\lambda)\; d\alpha(\lambda) = \sum_{i=1}^{N} w_i f(t_i) + \sum_{j=1}^{M} v_j f(z_j) + R[f]$$

$$where \qquad \text{Eq. (2)}$$

$$R[f] = \frac{f(\eta)^{(2N+M)}}{(2N+M)!} \int_a^b \prod_{j=1}^{M} (\lambda - z_j) [\prod_{i=1}^{N} (\lambda - t_i)]^2 \, d\alpha(\lambda), a < \eta < b$$

The above approximations of the integral when $M=0$ is called the Gauss formula, and when $M=1,2$ Gauss-Radau and Gauss-Lobatto respectively. When $M=1$, $z_1=a$ or $z_1=b$ (for upper or lower bounds respectively) and when $M=2$, $z_1=a$, $z_2=b$ (for an upper bound). The nodes $t$ and $z$ can be obtained by a polynomial decomposition of the integral in terms of $p_i(\lambda)$; moreover a set of orthogonal polynomials provides a 3-term recursion relationship for easy evaluations. This means the recurrence coefficients can be represented in a matrix of symmetric tri-diagonal form as in Eq. (3); the crucial observation being that these can be trivially extracted from the CG iterates, resulting in negligible addition of computation over the iterative solution. In more generality, the CG algorithm can be described as a minimization of the polynomial relation $\|x - x_k\|_A = \min_{p_k} \|p_k(A)(x - x_o)\|_A$.

$$\text{Given} \int_a^b p_i(\lambda)\ p_j(\lambda) d\alpha(\lambda) = \begin{cases} 0 \text{ when } i \neq j \\ 1 \text{ when } i = j \end{cases}$$

$\gamma_i p_i(\lambda) = (\lambda - \omega_i) p_{i-1}(\lambda) + \gamma_{i-1} p_{i-2}(\lambda), \quad i = 1,2,...N$ when normalized such that

$\int d\alpha = 1, p_0(\lambda) = 1, p_{-1}(\lambda) = 0.$

$$\Rightarrow \lambda P_{N-1}(\lambda) = C_N P_{N-1}(\lambda) + \gamma_N p_N(\lambda) e_N, \qquad \text{Eq. (3)}$$

where $e_N^T = [0,0,......1], P_{N-1}(\lambda)^T = [p_1(\lambda)\ p_2(\lambda) ... p_{N-1}(\lambda)]$

$$C_N = \begin{bmatrix} \omega_1 & \gamma_1 & 0 & 0 \\ \gamma_1 & \omega_2 & . & 0 \\ 0 & . & . & \gamma_{N-1} \\ 0 & 0 & \gamma_{N-1} & \omega_N \end{bmatrix}$$

Thus the eigenvalues of $C$ are the nodes of the Gauss formula, and the $C_R$, $C_L$ (for Radau and Lobatto formulas) can be computed from $C$ as shown elsewhere [5, 9, 10]. The above coefficients can be computed from the CG coefficients $\alpha_k$ and $\beta_k$ as below.

$$\omega_k = \frac{1}{\alpha_{k-1}} + \frac{\beta_{k-1}}{\alpha_{k-2}} \text{ and } \gamma_k = \frac{\sqrt{\beta_k}}{\alpha_{k-1}} \text{ where } \beta_0 = 0, \alpha_{-1} = 1.$$

where

$$\alpha_{k-1} = \frac{r_{k-1}{}^T A r_{k-1}}{d_{k-1}{}^T A d_{k-1}}, x_k = x_{k-1} + \alpha_{k-1} d_{k-1} \qquad \text{Eq. (4)}$$

$$r_k = r_{k-1} - \alpha_{k-1} A d_{k-1}, \beta_k = \frac{r_k{}^T r_k}{r_{k-1}{}^T r_{k-1}}, d_k = r_k + \beta_k d_{k-1}$$

It is efficient to compute the $A$-norm of the error using an approximation for $C_N^{-1}(1,1)$ as shown before [9].

$$\|x - x_k\|_A^2 = \|r_o\|^2 [C_N^{-1}(1,1) - C_k^{-1}(1,1)] \qquad \text{Eq. (5)}$$

This algorithm [5] depends on the estimates of the extremal eigenvalues, and the discussed results have been limited to absolute errors and a priori known extremal eigenvalues ($a = \lambda_{min}, b = \lambda_{max}$). We find that the upper bounds estimated to be moderately sensitive to approximations of '$a$' (shown in figure (2)); the lower bounds and the Gauss rule estimates are not useful as stopping criterion in generality, and so are not discussed further in this work. In fact, even a small overestimation of the smallest eigenvalue can make the upper bounds lower than the actual error (or result in complex valued estimated norms), thus failing, while any large underestimation makes the bounds too conservative to be helpful. Though figure (2) shows the sensitivity of the upper bound for a vector $b$ given as a sum of canonical basis for generality, the sensitivity to the estimate of smallest eigenvalue in specific cases can be far greater than shown there.

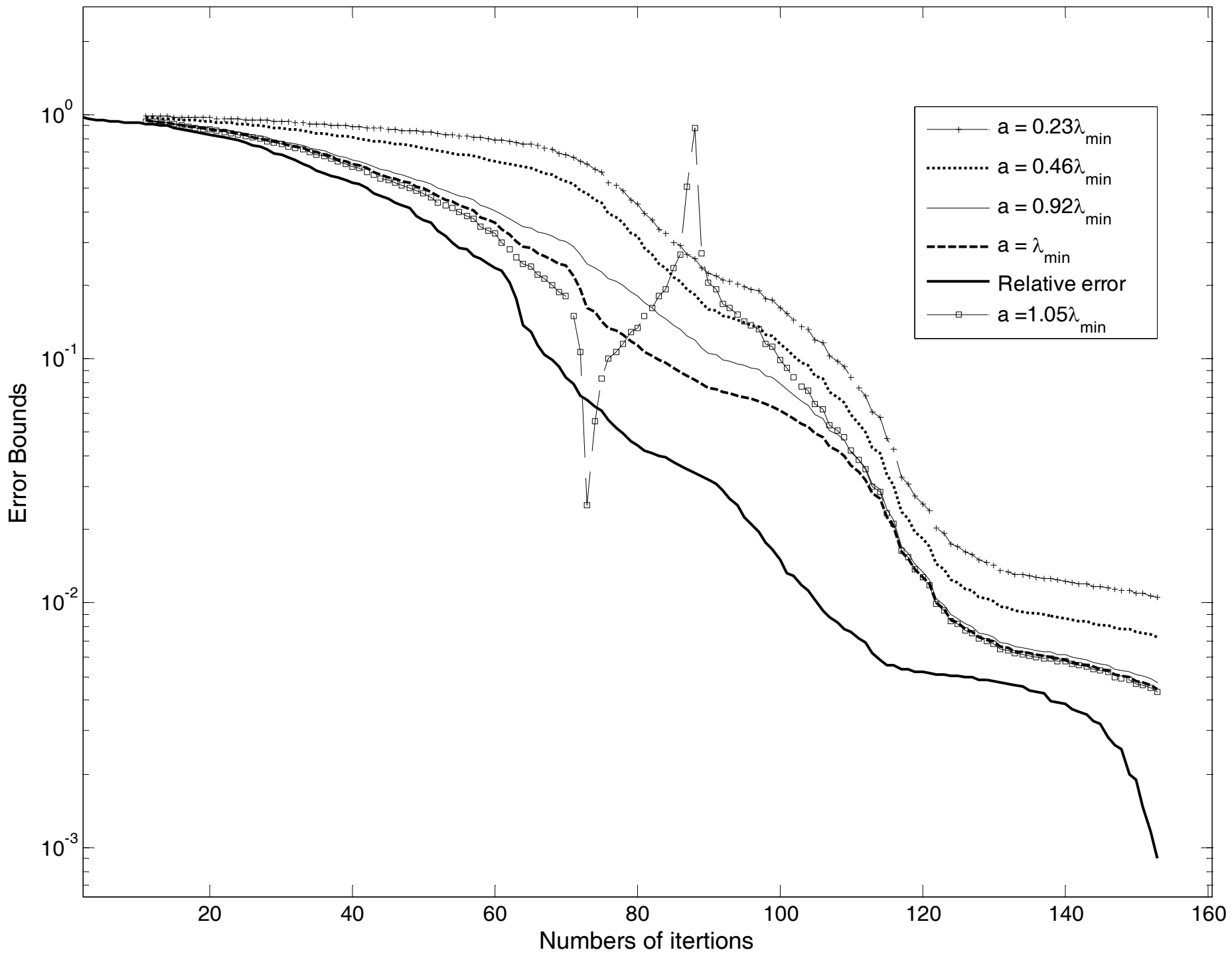

Figure 2: Sensitivity of the Gauss-Lobatto upper bounds of error (in $A$-norm) on the smallest eigenvalue '$a$' of SPD matrix $A$ (bcsstk05; math.nist.gov/MatrixMarket). Failure of the estimated bound is apparent when $a > \lambda_{min}$. Vector $b$ is a sum of canonical basis (all 1s) to maintain generality and $\kappa_A = 14281$, $N=153$.

## 3. Methods

We begin with showing that the eigenvalues of the tri-diagonal matrix $J_k$ constructed using CG iterates converge to those of $A$ as $k \rightarrow N$. Later we show that non-stationary convergence estimations using the smallest eigenvalue of $J_k$ can lower bind the smallest eigenvalue $a$ closely even if $k << N$; in fact with negligible computation; moreover, $J_k$ can be stored efficiently as a sparse matrix. One could include inverse iterations for the tri-diagonal $J_k$ at every CG iteration, nevertheless an approximation of '$a$' even in the same order of magnitude may mean one has to wait until $k \sim N$ for many matrices with large condition numbers. We suggest efficient non-stationary approximations of $a$, $b$ using the extremal eigenvalues of $J_k$ even as $k << N$. We also extend this method to estimation of the $l_2$ norm of the errors and the errors in relative terms. It is known that all symmetric matrices have a tri-diagonal symmetric form and this can be shown using orthogonal transformations.

***Theorem***: Any symmetric matrix $A$ has an orthogonal transformation to a tri-diagonal symmetric matrix $J$

Proof: Let $Q_2^T a_1$ be the orthogonal projection of $a_1$ the first column of $A$ such that $a_1^T Q_2 e_i = 0$ if $i > 2$ where $e_i$ denote the canonical basis vectors. Here orthogonality of $Q_2$ implies $a_i^T Q_2 Q_2^T a_j = a_i^T a_j$ *for* $i, j=1 \ldots N$, where $N$ is the dimension of $A$. $Q_2$ can be readily obtained by a Gram-Shmidt orthogonalization of $a_1$ with all $e_i$ for $i > 2$, or a Householder's reflection of $a_1$ over $e_1$, $e_2$.

Let $a_2$ be the second column of the symmetric matrix ($Q_2^T A\ Q_2$) and $Q_3^T a_2$ be the orthogonal projection of $a_2$ such that $a_2^T Q_3 e_i = 0$ *if* $i > 3$ $(i=1 \ldots.N)$. Continuing further in the same fashion to the orthogonal projection $Q_{N-1}^T a_{N-2}$ any symmetric matrix $A$ can be transformed into a symmetric tri-diagonal matrix $J = Q^T A Q = Q_{N-1}^T \ldots Q_3^T Q_2^T A\ Q_2 Q_3 \ldots\ Q_{N-1}$.

The above can also be shown using orthogonal transformations of the matrix $A$. Let $A^k = Q_k A^{k-1} Q_k^T$ be an orthogonal transformation such that $A^k_j$, the $j^{th}$ column of $A^k$ satisfies $e_i^T A^k_j = 0$ for all $j \leq k$ and all $i > j+1$. The resulting $A^{n-1} = J$ is the tri-diagonal matrix when $A^0 = A$ is the given symmetric matrix and $Q_k$ are the Householder's reflectors where

$$Q_k = \begin{bmatrix} I_k & 0 \\ 0 & H_{n-k} \end{bmatrix} \text{ and } H_{n-k} = I_{n-k} - \frac{2u_k u_k^T}{u_k^T u_k}, u_k = \left\| A_{k'}^{k-1} \right\| e_1 - A_{k'}^{k-1}, \qquad \text{Eq. (6)}$$

where $A_{k'}^{k-1} = [A_{k+1,k}^{k-1}, A_{k+2,k}^{k-1}, \ldots A_{n,k}^{k-1}]$ and note $H^2 = I$.

***Corollary***: Given $(J-\mu_i I)q_i = (Q^T A Q-\mu_i I)q_i = 0$ and $(A-\lambda_i I)p_i = 0$ where $\mu_i$, $q_i$ are eigenvalues and vectors of $J$; $\lambda_i$, $p_i$ are eigenvalues and vectors of $A$. Then it is evident that $\mu_i = \lambda_i$ and also, $q_i = Q^T p_i$ when $\lambda_i \neq 0$.

***Remark***: In more generality for a matrix $A$, if $A^T A = AA^T$, an eigenvalue revealing diagonal transformation $A = Q\Lambda Q^T$ exists and that is not proved here. The above orthogonal transformation of matrix $A$ into a tri-diagonal form is typically achieved by the Lanczos algorithm; it has a three-term recursion representation (like orthogonal polynomials) resulting in easier computation of a matrix $V_k$ including the first $k$ Lanczos vectors and transformed matrix $J_k$.

Moreover, Lanczos description helps in highlighting the equivalence of $C_k$ derived from the CG iterates as in Eq. (3), and the orthogonal transformations $A^k$ as below.

$$\text{Given } A^k = \begin{bmatrix} J_k & 0 \\ 0 & A^k_{n-k} \end{bmatrix} \text{ and } J_k = \begin{bmatrix} \omega_1 & \gamma_1 & 0 & 0 \\ \gamma_1 & \omega_2 & . & 0 \\ 0 & . & . & \gamma_{k-1} \\ 0 & 0 & \gamma_{k-1} & \omega_k \end{bmatrix}$$

$$V_k^T A V_k = J_k \text{ and } A V_k = V_k J_k + \gamma_{k+1} v_{k+1} e_k^T \qquad \text{Eq. (7)}$$

$$\Rightarrow A V_k = V_{k+1} J_k' \text{ where } J_k' = \begin{bmatrix} J_k \\ \gamma_{k+1} e_k^T \end{bmatrix}$$

Also known

$$\lambda P_k(\lambda) = P_{k+1}(\lambda) J_k' \text{ where } J_k' = \begin{bmatrix} J_k \\ \gamma_{k+1} e_k^T \end{bmatrix}$$

where

$P_k(\lambda)^T = [p_1(\lambda)\ p_2(\lambda) ... p_k(\lambda)]$ are a set of orthogonal polynomials such that Eq. (8)

$$\omega_k = \frac{1}{\alpha_{k-1}} + \frac{\beta_{k-1}}{\alpha_{k-2}} \text{ and } \gamma_k = \frac{\sqrt{\beta_k}}{\alpha_{k-1}} \text{ where } \beta_0 = 0, \alpha_{-1} = 1.$$

and $\alpha$,$\beta$ are CG iterates as in Eq.(4)

Thus when we use CG, we have all the elements of $J_k$, a partial orthogonal transformation of $A$ at the $k^{th}$ iteration. The *largest* eigenvalue of $J_k$ is known to converge to that of $A$ quickly, but it is the *smallest* eigenvalue that may need almost $N$ iterations to converge completely. The largest eigenvalue of $J_k$ is a very good approximation of the required largest eigenvalue $b$ even as $k<<N$, and so we will not discuss it further. In general, most iterative solutions of linear systems can be considered effective only when the number of iterations taken to solve a large problem is much less than $N$. We compute the extremal eigenvalues of $J_k$ and subsequently underestimate the smallest eigenvalue of $A$ even as $k<<N$ using a non-stationary function, with the knowledge that the extremal eigenvalues of $J$ converge to that of $A$ with increasing $k$. We can evaluate the largest and smallest eigenvalues of $J_k$ in < $O(k^2)$ operations using methods (available in LAPACK libraries) that are not discussed here.

If $f_k$ is the smallest eigenvalue of $J_k$, it is guaranteed that $0 < f_{k+1} \leq f_k$ where $f_{k+1}$ is the smallest eigenvalue of $J_{k+1}$, given the matrix $A$ is positive definite. Any other correlation between $f_k$ and $f_{k+1}$ is not known to satisfy the generality of spectral properties required of any symmetric positive definite matrix $A$. But it is observed for all matrices studied that the finite difference of the sequence of the smallest eigenvalues of $J_k$ is in fact *nominally* decrease with $k$ (as defined by $\beta < 1$ in the statement below). Moreover, any sequence with a nominally decreasing finite difference can be represented by a sequence with piecewise strictly decreasing finite differences. Hence for sequences with decreasing finite differences, we can find an appropriate non-stationary function with some rigor that is computationally trivial to update. We need a look-ahead *underestimate* of smallest eigenvalue of $A$ that modulates $f_k$ by a factor between 1 and 0, depending on the convergence rate of $f_k$. The exponential function (with exponent ≤ 0) is then a very suitable look-ahead to lower bind the smallest eigenvalue of $A$ using the eigenvalue convergence of $J_k$. (Note: In the case of a continuous function or a discrete sequence, unless $f(s)$ or $f(z)$ the Laplace transform of the function or the Z-transform of the sequence can be well approximated in the limits $s \to 0, z \to 1$, the final value theorem is not useful).

***Statement***: Let $f$ be any real valued sequence of numbers $\forall f_{k-1} \geq f_k > 0$ for $k = 1 \ldots N$ and $\frac{1}{N-1}\sum \frac{f_k - f_{k+1}}{f_{k-1} - f_k} = \beta < 1$ (defined for all $f_{k-1}$-$f_k \neq 0$ in the sequence). '$f_k$' is a sequence of positive real numbers for integers $k = 1 \ldots N$ such that both its value and the *nominal* finite difference of the sequence are not increasing with $k$. Given $f_N = a$, find the positive valued function $g_k$ in $f_k e^{-g_k(f_i,k,N)} = \tilde{a}_k$ where $i = 1 \ldots k < N$ such that

(1) $\tilde{a}_k \leq a$ and (2) $\frac{\sum_{k=1}^{N}(a - \tilde{a}_k)^2}{\sum_{k=1}^{N}(f_k - a)^2} = \delta << 1$

For any arbitrary sequence $f_k$, simultaneous enforcement of both conditions (1, 2) can not be guaranteed for all $k$ strictly; then minimizing $\delta$ in (2) while satisfying (1) for almost all $k$ is a reasonable objective for the required lower bound estimate of $f_N = a$.

***Proposition***:

Let $f_k \geq f_{k+1} > 0$ and $\frac{f_k - f_{k+1}}{f_{k-1} - f_k} \leq 1$ for $k = 1 \ldots N$.

Then a function $f_k e^{-g_k(f_i,k,N)} = f_k \exp(-\alpha_k \times [\frac{N-k}{N}]) = \tilde{a}_k$ ensures there exists $\alpha_k$ such that

$$\frac{\sum_{k=1}^{N}(a-\tilde{a}_k)^2}{\sum_{k=1}^{N}(f_k - a)^2} = \delta << 1 \text{ when } f_N = a.$$

*Proof*: We strictly enforce $f_N$ on the final value of the sequence approximated $\tilde{a}_k = f_k e^{-g_k(f_i,k,N)}$ for any sequence $f$ such that $\tilde{a}_N = f_N = a$ for $k=N$. This implies $g_N = 0$ and this is easily done for all $f$ when $g_k(f_i,k,N) = \alpha_k(f_i) \times F(k,N)$ where $i=1...k$, a possible separation of variables for function $g$, and $F(k,N) \approx d \times [N-k]$ to the first order. The choice for the multiplicative factor $d$ can either be a constant such as $N^{-1}$, a variable such as $k^{-1}$, or even a random variable. If $F(k,N) = [N-k]/k$, $\tilde{a}_k$ is a very sensitive to $\alpha_k$ when $k << N$, but becomes stiff as $k/N \to 1$, and on the other hand if $F(k,N) = [N-k]/N$ it is relatively stiff throughout; we present results of both cases. Then it is evident that when $\tilde{a}_k = f_k \exp(-\alpha_k \times [N-k]/N)$, for any $f$, $\alpha_k = \frac{N}{N-k}\ln(\frac{f_k}{f_N})$ and $\tilde{a}_k = a \Rightarrow \delta = 0$.

But to find $\alpha_k$ given only $f_i, i \le k < N$, we can minimize the error between the non-stationary estimates of the current value using the previous values of the sequence, given by $(f_i \exp(-\alpha_k \times [k-i]/k) - f_k)^2$. When a strictly decreasing finite difference of the sequence is not guaranteed and it has decreasing finite differences only piecewise, simultaneous enforcement of both conditions 1, 2 of statement (1) are not possible for all $k$ strictly. This is because then any estimate of $\alpha_k$ will not monotonically reduce with $k$, and a bound on $\delta$ is difficult for such an arbitrary sequence. In such cases, an estimate of $\tilde{a}_k$ can be highly oscillatory because of the large $N$-$k$ that is to be *looked ahead* when $1<k<<N$. The logarithmic scaling in the regression, (i.e.) the approximations of $\ln(f_i / f_k)$ for all $i$ such that $i+m>k$, can not only reduce large such variations in $\tilde{a}_k$ but result in a much cheaper computation as Eq. (9) in an estimation $\underline{\alpha}_k$ for the underestimate of $f_N = a$. The computational cost of the regression in evaluating $\underline{\alpha}_k$ is only $O(m)$ where the last $m$ points in the sequence are used to minimize the error in exponential fits of sequence $f$. This helps in limiting the total computation of the proposed error estimators to much less than an iteration of hosting CG algorithm which is $O(N^2)$, and thus a recursive estimate of extremal eigenvalues can be evaluated in all iterations of CG. The following is then a trivial update of the estimate of the smallest eigenvalue of matrix *A*.

$$\frac{d\sum_{i=k-m}^{k-1}(\ln(f_i/f_k)-\underline{\alpha}_k\times[k-i]/k)^2}{d\underline{\alpha}_k}=0$$

$$\Rightarrow \underline{\alpha}_k=\frac{\sum_{i=k-m}^{k-1}\frac{k-i}{k}\times\ln\frac{f_i}{f_k}}{\sum_{i=k-m}^{k-1}(\frac{k-i}{k})^2} \text{ and } \tilde{a}_k=f_k\exp(-\underline{\alpha}_k\times[\frac{N-k}{N}]) \qquad \text{Eq. (9)}$$

Though Figure (3) plotted in logarithmic scale might mislead one to conclude that these look-ahead underestimates are useful only for the initial iterations, it should be reminded that even a small but continued overestimate of smallest eigenvalue '*a'* can fail the upper bounds on error (as shown in figure (2)). Hence these exponential lower bounds of the $\lambda_{min}$ are useful for tight and assured upper bounds of relative error even as $k \rightarrow N$ by ensuring a tight underestimate of the smallest eigenvalue of *A*. Also, iterative solutions of linear systems are considered effective for a problem only when a stopping $k/N < 0.1$, when an underestimate of $\lambda_{min}$ is critical for the estimate of upper bound of errors.

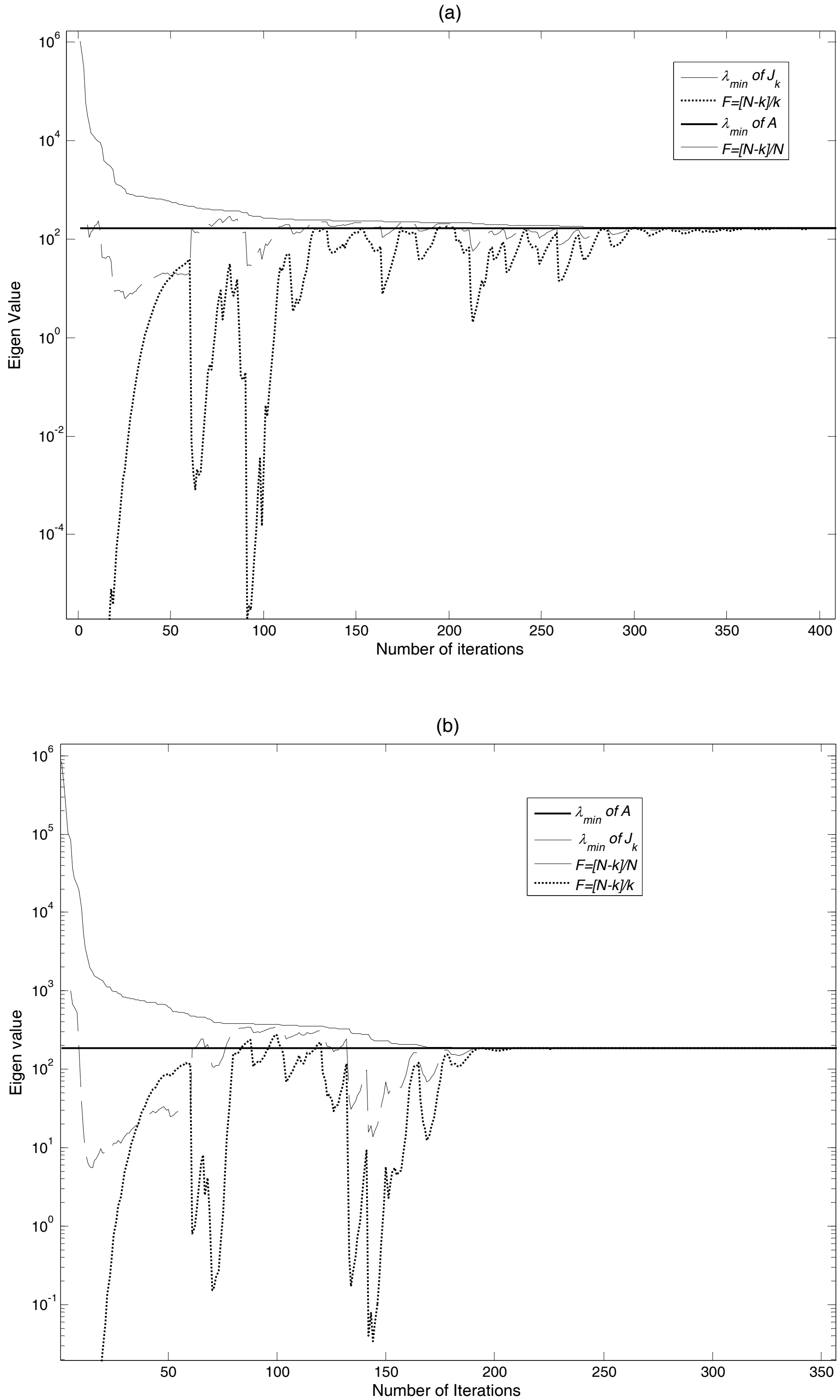

Figure 3: $\lambda_{min}$ estimates using the proposed extremal eigenvalue approximation for SPD matrices (math.nist.gov/MatrixMarket) along with CG: (a) bcsstm19; $\kappa_A$ = *233734*, *N=817* and (b) bcsstm20; $\kappa_A$ = *255380*, *N=485* when vector *b* is a sum of canonical basis for generality and *m = 10*.

## 4. Numerical results and Discussion

The estimated extremal eigenvalues are then introduced into the quadrature based estimates for both relative errors in $l_2$ norm and $A$-norm. As mentioned before, the $l_2$ norm estimates are not guaranteed to be upper bounds but nevertheless are close enough to the actual relative errors to be used as a stopping criterion. Actual errors are evaluated using the final $x$ computed after $N$ iterations and $x_k$ at $k^{th}$ iteration. This is valid when the relative residues iterate to $\|r\|/\|b\| << 1$, indicating convergence of the CG algorithm. Once the $A$-norm of the absolute error is estimated, the $A$-orthogonality of the CG iterates can be invoked to estimate the errors in relative terms. $||x||_A^2 = ||x - x_k + x_k||_A^2$ and thus $||x||_A = (||x - x_k||_A^2 + ||x_k||_A^2)^{1/2}$ when $<x-x_k, Ax_k> = 0$; hence $||x - x_k||_A^2 / ||x||_A^2$ can be well approximated.

The evaluation of bounds of $l_2$ norm of the error, $||x-x_k||^2 = r^T A^{-2} r$ using the quadrature rules involves additional difficulties; especially in ascertaining if the evaluations are upper or the lower bounds [11]. The $l_2$ norm of the absolute error in terms of the matrix $J$ is given by Eq. (10) and approximations of the first term is especially non-trivial. Unlike there where Gauss rules were used, we have used the extremal eigenvalue approximations in estimating Gauss-Radau lower bounds of $\|\varepsilon\|_A^2$; this used in Eq. (10) ensures we have an upper bound of $l_2$ norm most likely. The estimates of relative error in these cases is well approximated by $||x-x_k||/|| x|| \approx ||x-x_k||/|| x_k||$ when it is $<< 1$, in spite of oscillations in $|| x_k||$, and we show a few examples in figure (5) that use the extremal eigenvalue approximations in the computation.

$$\|\varepsilon\|^2 = \|r_o\|^2 [e_1^T J_N^{-2} e_1] - [e_1^T J_k^{-2} e_1] - 2\frac{[e_k^T J_k^{-2} e_1]}{[e_k^T J_k^{-1} e_1]}\|\varepsilon\|_A^2 \qquad \text{Eq. (10)}$$

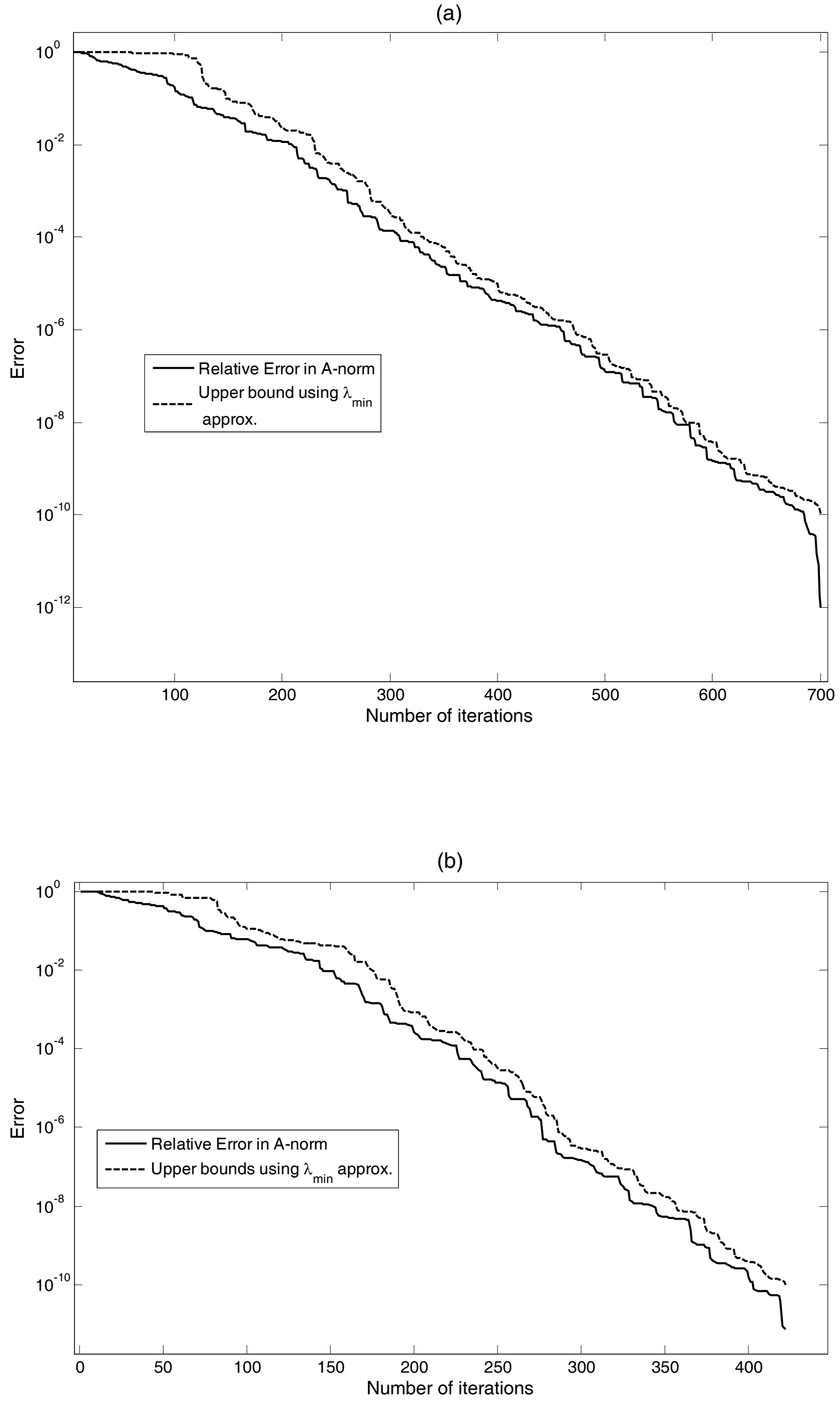

Figure 4: Gauss-Lobatto upper bounds of relative errors in *A*-norm using the proposed extremal eigenvalue approximation for SPD matrices: (a) bcsstm19; $\kappa_A$ = *233734*, *N=817* and (b) bcsstm20; $\kappa_A$ = *255380*, *N=485* (math.nist.gov/MatrixMarket) when vector *b* is a sum of canonical basis for generality.

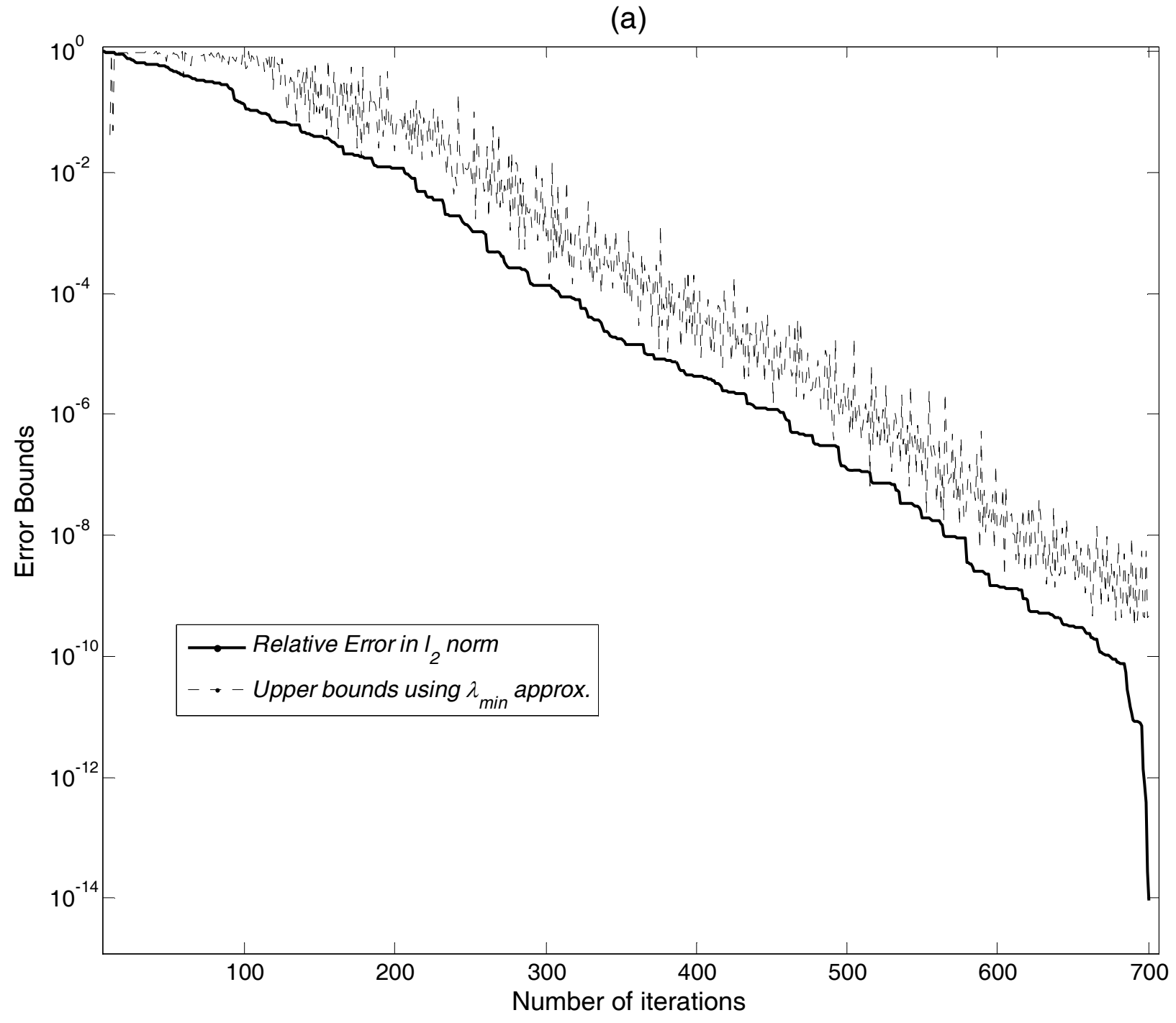

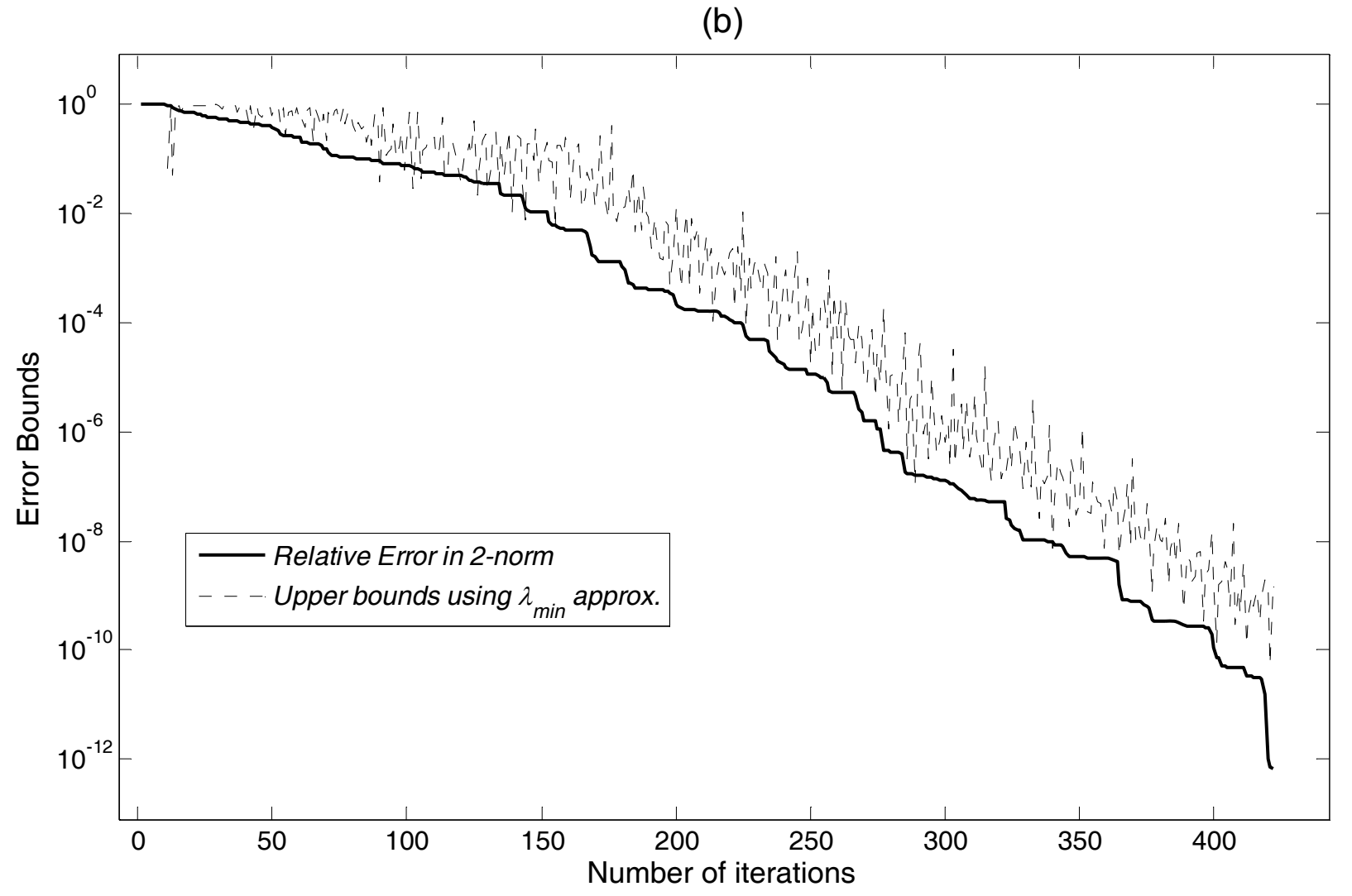

Figure 5: Estimates of $l_2$ norm of relative error using the proposed extremal eigenvalue approximation and the Gauss-Radau rules for SPD matrices: (a) bcsstm19; $\kappa_A = 233734$, $N=817$ and (b) bcsstm20; $\kappa_A = 255380$, $N=485$ (math.nist.gov/MatrixMarket) when vector $b$ is a sum of canonical basis for generality.

## 5. Conclusions

The significance of error estimators for even problems with a moderately high condition number ($>10^2$) is evident and was emphasized by a few general examples in the introduction. A stopping criterion by the residue can be orders of magnitude ($\sim\kappa$) conservative or optimistic and thus computationally either very costly or error prone for large problems. Though results of error estimation without a priori knowledge of eigenvalues have been shown here for SPD matrices and the CG algorithm, other iterative algorithms for generalized matrices can as well include such eigenvalue approximations and error estimations. There can be matrices of large condition numbers and clustered eigenvalues where Lanczos algorithm may need more than $N$ iterations due to finite numerical precision (in converging to an tri-diagonal matrix that includes *all* eigenvalues of $A$); the methods presented here are not in anyway limited by such matrices and the number of iterations of the CG algorithm.

**Acknowledgements**

DAS thanks IISc Mathematics Initiative (IMI) for the generous support of its graduate students. Both DAS and MV thank Prof. Srinivas Talabatulla for his mentorship of DAS.